\documentclass[12pt]{amsart}
\usepackage[english]{babel}
\usepackage{amssymb}

\newtheorem{theorem}{Theorem}
\newtheorem{conjecture}{Conjecture}

\newcommand{\Z}{\mathbb Z}
\newcommand{\Q}{\mathbb Q}
\newcommand{\R}{\mathbb R}
\newcommand{\Co}{\mathbb C}
\DeclareMathOperator{\Hom}{Hom}
\DeclareMathOperator{\E}{E}
\DeclareMathOperator{\ord}{ord}
\DeclareMathOperator{\Aut}{Aut}

\DeclareMathOperator{\To}{\longrightarrow}

\title[Generalization of the ``Stark unit'']{Generalization of the ``Stark unit'' for abelian L-functions with multiple zeros}
\author{M.~Vlasenko}
\begin{document}
\maketitle

Let $K / k$ be a normal extension of number fields with abelian Galois group $G$. Suppose there is an infinite place $v$ of $k$ such that the set of nontrival characters $\chi : G \To \Co^{\times}$ which are ramified at all infinite places except $v$ is nonempty. Then the corresponding Artin L-functions $L(\chi, s)$ have zeros of the first order at $s=0$, and if we fix any embedding of the bigger field $K$ to $\Co$ which extends $v$ then Stark's conjecture predicts existence of a ``unit'' $\varepsilon \in \Q \otimes_{\Z} O_K^{\times}$ such that for every such $\chi$
\[
L'(\chi, 0) = \sum_{g \in G} \chi(g) \log |g \varepsilon|.
\]
We refer to the book~\cite{Tate} concerning the formulation and consistency of the conjecture. It is stated there for a character of an arbitrary Galois representation, but the case of abelian character would imply the general one. 

In this note we fix a nonnegative integer $r \ge 0$ and consider those $\chi$ for which the L-function has zero of order exactly $r$ at $s=0$. Let
\[
L(\chi, s) = c(\chi) s^r + o(s^r), \;\;\; s \rightarrow 0
\]
with $c(\chi) \ne 0$. First we show that the leading Laurent coefficients $c(\chi)$ for such $\chi$ can be all described by means of a single element $\varepsilon_r \in \R \bigwedge_{\Z[G]}^r O_K^{\times}$ of the $\R$-span of the $r$-th exterior power of the module of units $O_K^{\times}$ taken in the category of $\Z[G]$-modules. Further, we prove that the Stark conjecture is equivalent to the statement that
\[
\varepsilon_r \in \Q \bigwedge_{\Z[G]}^r O_K^{\times}.
\]
Our main result is~Theorem~2. Following~\cite{Tate}, we consider here slightly more general L-functions $L_S(\chi, s)$ where $S$ is a set of places of $k$ containing all the infinite places $S_{\infty} \subseteq S$. If one puts $S=S_{\infty}$ then these are the ordinary Artin L-functions as above, namely $L(\chi,s)=L_{S_{\infty}}(\chi, s)$. In general $\varepsilon_r$ belongs to $\R \bigwedge_{\Z[G]}^r U$, where $U \subset K^{\times}$ is the module of $S$-units of the bigger field $K$. If $S=S_{\infty}$ then $U = O_K^{\times}$. 

We should mention that such elements already appeared in~\cite{Rubin}, but under very special assumptions on the set of places $S$  (Hypotheses 2.1 of~\cite{Rubin}). Although these assumptions allow to refine the Stark conjecture, we find it convenient from the computational point of view to keep the set $S$ as small as possible. Indeed, when we increase $S$ then the order of zero of $L_S(\chi, s)$ at $s=0$ either stays constant or increases. The module $U$ of S-units also grows. Hence if one really wants to verify the Stark conjecture numerically, the most simple choice will be $S=S_{\infty}$. We consider one numerical example below. In view of the above-said, the present note may be considered as a naive version of~\cite{Rubin}, where we do everything what can be done simply by linear algebra.

\section{The conjecture}
We start with a reformulation of the abelian Stark conjecture. Let $K / k$ be a normal extension of number fields with an abelian Galois group $G$, $S$ be a set of places of $k$ containing all the infinite places $S_{\infty} \subset S$ and $S_K$ be the set of all places of $K$ lying above places in $S$. Let $Y$ be the free abelian group on $S_K$ and
\[
X = \left\{ \sum_{w \in S_K} x_w w \; \Big| \; x_w \in \Z, \;\sum x_w = 0 \right\} \subset Y.
\]
Let 
\[
U = \left\{ x \in K^{\times} \; \Big| \; \ord_w (x)=0 \;\;\; \forall w \notin S_K \right\}
\]
be the group of $S$-units in $K$. We have a $G$-module isomorphism 
\[
\log: \R U \To \R X
\]
\[
\log(x) = \sum_{w \in S_K} \log|x|_w \; w.
\] 
Also $\Q X \cong \Q U$ as $G$-modules (see~\cite{Tate}, p.26), but there is no canonical choise for this isomorphism. 

The modified Artin L-function attached to a character $\chi: G \To \Co^{\times}$ is defined for $\Re(s) > 1$ as
\begin{equation}\label{Artin_L}
L_S(\chi, s) = \prod_{v \notin S, \chi|_{I_v} \equiv 1} \left( 1-\chi(\sigma_v\right) Nv^{-s})^{-1},
\end{equation}
where for every finite place $v$ of $k$ we denote by $G_v, I_v$ and $\sigma_v\in G_v /I_v$ the decomposition group, the inertia group and the Frobenius element correspondingly. The order of zero at $s=0$ of this L-function is 
\[
r(\chi) = \# \left\{ v \in S \;\Big| \; \chi|_{G_v} \equiv 1 \right\}
\]
for $\chi \ne 1$ and $r({\bf 1})=\#S - 1$ (Proposition II.3.4 in \cite{Tate}). Here ${\bf 1}$ is the trivial character. Let $c(\chi)$ be the leading Laurent coefficient of the L-function~(\ref{Artin_L}) at $s=0$:
\[
L_S(\chi,s)=c(\chi) s^{r(\chi)}+o(s^{r(\chi)}), \;\;\; s \rightarrow 0.
\]

Let us fix a nonnegative integer $r$ such that there are characters $\chi$ with $r(\chi)=r$. In particular, should be $r \le \#S$. Let $C$ be a nonempty set of characters satisfying the two conditions:
\begin{itemize}
\item[(i)] if $\chi \in C$ then $\chi^{\alpha} = \alpha \circ \chi \in C$ for every $\alpha \in Aut(\Co/\Q)$; 
\item[(ii)] $r(\chi) = r$ for all $\chi \in C$.
\end{itemize}
For any $\chi$ the element $\E_{\chi} = \frac 1{\# G} \sum_{g \in G} \overline \chi(g) g$ is an idempotent in $\Co[G]$. Let 
\[
\E_C = \sum_{\chi \in C} \E_{\chi}.
\]
It is an idempotent since $\E_{\chi} \E_{\chi'} = 0$ when $\chi \ne \chi'$, and it belongs to $\Q[G]$ due to the property~(i) of the set $C$. Consider also the generalized ``Stickelberger element'' 
\[
\Theta_C = \sum_{\chi \in C} c(\overline \chi) \E_{\chi}.
\]
It obviously belongs to the subspace $\E_C \Co[G]$, but since $c(\overline \chi)=\overline{c(\chi)}$ we have $\Theta_C \in \E_C \R[G]$ in fact. 

Let $\log^{(r)} : \R \bigwedge^r U \To \R \bigwedge^r X$ be the $G$-module isomorphism induced by $\log : \R U \To \R X$, where exterior products are always taken in the category of modules over the commutative ring $\Z[G]$. If $r=0$ this means $\bigwedge^0 U = \bigwedge^0 X = \Z[G]$ and $\log^{(0)} = 1$.

\begin{conjecture} \footnote{If $r=0$ the conjecture states $\Theta_C \in E_C \Q[G]$. The same claims the Stark conjecture for all $\chi \in C$, which is true (Theorem III.1.2 of \cite{Tate}). We assume $r>0$ further in this note.}
\label{Stark_over_Q} Let $C$ be a nonempty set of characters satisfying the conditions~(i) and~(ii) from above. Then 
\[
\Theta_C \; \Q \bigwedge^r X = \Q \log^{(r)} \left( \E_C \bigwedge^r U \right).
\] 
\end{conjecture} 
Notice that both $\Theta_C \; \bigwedge^r X$ and $\log^{(r)} \left( \E_C \bigwedge^r U \right)$ are maximal discrete lattices in the real vector space $\E_C \R \bigwedge^r X$ of dimension $\# C$. Indeed, $\E_C \bigwedge^r X$ is obviously a maximal lattice in $\E_C \R \bigwedge^r X$, and $\Theta_C$ is an invertible transformation of this space since $\det(\Theta_C) = \prod_{\chi \in C} c(\chi) \ne 0$. Analogously, $\E_C \bigwedge^r U$ modulo torsion is maximal in $\E_C \R \bigwedge^r U$, and $\log^{(r)}: \E_C \R \bigwedge^r U \To \E_C \R \bigwedge^r X$ is invertible. 

The conjecture claims these lattices are commensurable.

\begin{theorem}\label{equivalence}
Conjecture~\ref{Stark_over_Q} is equivalent to the Stark conjecture ``over $\Q$'' (Conjecture I.5.1 of \cite{Tate}) for all characters $\chi \in C$.
\end{theorem}

Let us choose any isomorphism $f: \Q X \To \Q U$. Then the Stark regulator for each $\chi$ is defined as $R(\chi, f) = \det(\log \circ f |_{\Hom_G(\overline \chi, \Co X)})$. Put $A(\chi,f)=\frac{c(\chi)}{R(\chi,f)}$. Stark's conjecture claims that $A(\chi^{\alpha},f)=A(\chi,f)^{\alpha}$ for every $\alpha \in \Aut(\Co/\Q)$. We can equivalently state it as 
\begin{equation}\label{Stark_rational}
A = \sum_{\chi \in C} A(\overline \chi, f) \E_{\chi} \in \Q[G].
\end{equation}
Let us introduce also $R = \sum_{\chi \in C} R(\overline \chi, f) \E_{\chi}$. Then $\Theta_C = A \; R$.

\begin{proof}[Proof of Theorem~1.] For any $x_1, \dots, x_r \in X$
\[
R(\overline \chi, f) \E_{\chi} x_1 \wedge \dots \wedge x_r = \E_{\chi} \log \circ f(x_1) \wedge \dots \wedge \log \circ f(x_r)
\]
by the definition of $R(\overline \chi, f)$. Indeed, $\Hom(\chi, \Co X) \cong \E_{\chi} \Co X$ and dimension of this space is $r$. Therefore
\begin{equation}\label{regulator_transform}
R \bigwedge^r X = \log^{(r)} \left( \E_C \bigwedge^r f(X) \right),
\end{equation}
and $\Q \log^{(r)} \left( \E_C \bigwedge^r U \right) = R \; \Q \bigwedge^r X$. Also $\Theta_C \Q \bigwedge^r X = R \left( A \; \Q \bigwedge^r X \right)$. Since multiplication by $R$ is an invertible transformation of $\E_C \R \bigwedge^r X$, Conjecture~\ref{Stark_over_Q} is now equivalent to the statement $\E_C \Q \bigwedge^r X = A \; \Q \bigwedge^r X$. This is in turn equivalent to $A \in \E_C \Q[G]$. Recall that we reduced the  Stark conjecture to the same statement~(\ref{Stark_rational}).
\end{proof}

We have already remarked that Conjecture~\ref{Stark_over_Q} claims commensurability of the two lattices $\Theta_C \bigwedge^r X$ and $\log^{(r)} \left( \E_C \bigwedge^r U \right)$. Let us show that the ratio of covolumes of these lattices is a rational number. Let $f: \Q X \To \Q U$ be the isomorphism from above. Then $\log^{(r)} \left( \E_C \bigwedge^r f(X) \right)$ is commensurable with $\log^{(r)} \left( \E_C \bigwedge^r U \right)$, and
\[
\Theta_C \bigwedge^r X = A R \bigwedge^r X = A \log^{(r)} \left( \E_C \bigwedge^r f(X) \right) 
\]
due to~(\ref{regulator_transform}). The ratio of covolumes of $\Theta_C \bigwedge^r X$ and $\log^{(r)} \left( \E_C \bigwedge^r f(X) \right)$ is therefore $\det(A) = \prod_{\chi \in C} A(\chi, f)$. 
For a character $\psi: G \To \Co$ of an arbitrary representation of $G$ one can define the Artin L-function $L_S(\psi, s)$ and introduce the numbers $c(\psi)$, $R(\psi, f)$ and $A(\psi, f) = \frac{c(\psi)}{R(\psi,f)}$. Since the latter numbers satisfy $A(\psi_1+\psi_2)=A(\psi_1)A(\psi_2)$, 
\[
\det(A) = \prod_{\chi \in C} A(\chi, f) = A \left( \sum_{\chi \in C} \chi \;, \; f \right).
\]
Due to the property~(i) of $C$ the character $\psi = \sum_{\chi \in C} \chi$ takes rational values on $G$. For such characters the Stark conjecture is known to be true (Corrolary II.7.4 in~\cite{Tate}).\footnote{We refer here to the truth of the Stark conjecture for a rational-valued character, although one can find simple arguments in this case. In fact our $\psi$ is a sum of permutation characters. A permutation character is the one induced from the trivial character on a subgroup. Since $A(\mathrm{Ind} \; \chi, f) = A(\chi, f)$, the Stark conjecture is true for $\psi$ because it is true for a trivial character.} It means $A(\psi,f) \in \Q$, hence our claim is proved.
 
\section{Generalization of the "Stark unit"}

Let $S=\{ v_1,\dots, v_n\}$ and choose some $w_i \in S_K$ above each $v_i \in S$. For $w \in S_K$ let $w^* \in \Hom_G(Y, \Z[G])$ be the map 
\[
w^* \left( \sum x_{w'} w' \right) = \frac 1{\#G_{v_i}} \sum_{g \in G} x_{g w} \; g.
\]
Then we can consider $w_{i_1}^* \wedge \dots \wedge w_{i_r}^* \in \Hom_G(\bigwedge^r_G Y, \Z[G])$ for any $i_1 < \dots < i_r$. This homomorphism maps $y_1 \wedge \dots \wedge y_r$ to $\det (w_{i_k}^*(y_l))_{k,l=1}^r$.
 
Let $\E_r = \underset{\chi: r(\chi)=r} \sum \E_{\chi}$. It is an idempotent in the group algebra $\Q[G]$. Recall that $U \subset K^{\times}$ is the group of $S$-units of $K$.

\begin{theorem}\footnote{Let us write for the moment $\varepsilon^K_r$, $w_i^K$ and $U_K$ to indicate the relation to $K/k$. Let us take some intermediate extension $k \subset F \subset K$ and choose for the places $w_i^F$ exactly the ones under $w_i^K$.  Then one can prove that $\varepsilon^F_r = \mathcal N^{(r)}(\varepsilon^K_r)$, where $\mathcal N^{(r)}: \R \bigwedge^r U_K \To \R \bigwedge^r U_F$ is the map induced by the norm $\mathcal N: U_K \To U_F$.} There is a unique element $\varepsilon_r \in \E_r \; \R \bigwedge^r U$ such that for every $\chi \ne {\bf 1}$ with $r(\chi)=r$
\[
c(\chi) = \overline{\chi} \left( w_{i_1}^* \wedge \dots \wedge w_{i_r}^* \left(\log^{(r)}(\varepsilon_r)\right) \right)
\]
where $\{ i_1 < \dots < i_r \} = \{ i : \chi|_{G_{v_i}} \equiv 1\}$, and (in case $r = n-1$)
\[
c({\bf 1}) =  {\bf 1} \left( (-1)^i w_1^* \wedge \dots \wedge w_{i-1}^* \wedge w_{i+1}^* \wedge \dots \wedge w_n^* \left(\log^{(r)}(\varepsilon_r)\right) \right)
\]
with an arbitrary $i$. 

Conjecture~\ref{Stark_over_Q} for the set of characters $C_r=\{\chi: r(\chi)=r\}$ is true if and only if $\varepsilon_r \in \E_r \; \Q \bigwedge^r U$.
\end{theorem}

\begin{proof} We want to construct a certain $G$-module isomorphism 
\[
\Omega : \E_r \Q \bigwedge^r X \; \overset{\sim} \To \; \E_r \Q [G].
\]
Suppose first $r \ne n-1$, so that all $\chi \in C_r$ are nontrivial. Let us check that
\[
\Omega = \underset{i_1 < \dots < i_r}  \sum w_{i_1}^* \wedge \dots \wedge w_{i_r}^*  \;\;\; \Big|_{\E_r \Q \bigwedge^r X}
\]
is such an isomorphism. Since its image obviously belongs to $\E_r \Q[G]$, it is enough to prove that $\Omega$ gives an invertible map from $\E_r \Co \bigwedge^r X$ to $\E_r \Co[G]$ or, equivalently, from each $\E_{\chi} \Co \bigwedge^r X$ to $\E_{\chi} \Co[G]$. The last space is of dimension 1. Notice that if $\chi|_{G_{v_i}}$ is nontrivial then $\E_{\chi} w_i = 0$ in $\Co Y$. If $\chi|_{G_{v_i}} \equiv 1$ then $\E_{\chi} w_i \in \E_{\chi} \Co X$. Let $\{ i_1 < \dots < i_r \} = \{ i : \chi|_{G_{v_i}} \equiv 1\}$. Then $\E_{\chi} w_{i_1}, \dots, \E_{\chi} w_{i_r}$ is a basis of $\E_{\chi} \Co X$, and therefore the space $\E_{\chi} \Co \bigwedge^r X = \Co \E_{\chi} w_{i_1} \wedge \dots \wedge w_{i_r}$ has dimension 1. Since $\Omega(\E_{\chi} w_{i_1} \wedge \dots \wedge w_{i_r}) = \E_{\chi}$ we are done.   

Let now $r=n-1$ and  
\[
\Omega = (\E_r-\E_{\bf 1})\underset{i_1 < \dots < i_r}  \sum w_{i_1}^* \wedge \dots \wedge w_{i_r}^* + \E_{\bf 1} (-1)^i \check w_i^* \;\;\; \Big|_{\E_r \Q \bigwedge^r X}
\]
where we denote $w_1^* \wedge \dots \wedge w_{i-1}^* \wedge w_{i+1}^* \wedge \dots \wedge w_n^*$ by $\check w_i^*$. Let us see that this map is independent of $i$ and it is an isomorphism from $\E_r \Q \bigwedge^r X$ to $\E_r \Q[G]$. We can tensor with $\Co$ again, and if $\chi \ne {\bf 1}$ then $\Omega$ is an isomorphism from each $\E_{\chi} \Co \bigwedge^r X$ to $\E_{\chi} \Co[G]$. $\E_{\bf 1} \Co \bigwedge^r X$ is generated over $\Co$ by $\E_{\bf 1} (w_2 - w_1)\wedge \dots \wedge (w_n - w_1) = \E_{\bf 1} \sum_j (-1)^{j-1} \check w_j$, where $\check w_j = w_1 \wedge \dots \wedge w_{j-1} \wedge w_{j+1} \wedge \dots \wedge w_n$. Since $\E_{\bf 1} \check{w}^*_i (\check w_j) = \E_{\bf 1}$ when $i=j$ and 0 otherwise, we are done.

With $\Omega$ constructed above we have an invertible map
\begin{equation}\label{omega_log_isom}
\Omega \circ \log^{(r)} : \E_r \R \bigwedge^r U \; \overset{\sim} \To \; \E_r \R[G].
\end{equation}
Notice that the property of $\varepsilon_r$ can be reformulated as 
\[
c(\chi)=\overline{\chi}(\Omega \circ \log^{(r)}(\varepsilon_r))
\]
for all $\chi$ with $r(\chi)=r$. Since $\Theta_r = \sum_{r(\chi)=r} c(\overline \chi) \E_{\chi}$ is a unique element in $\E_r \R [G]$ with the property $c(\chi)=\overline{\chi} (\Theta_r)$ for all such $\chi$, then $\varepsilon_r = (\Omega \circ \log^{(r)})^{-1} \Theta_r$ is a unique element satisfying the conditions.

Since $\Omega : \E_r \R \bigwedge^r X \To \E_r \R[G]$ is an isomorphism and it maps $\E_r \Q \bigwedge^r X$ to $\E_r \Q[G]$, it also maps $\Theta_r \Q \bigwedge^r X$ to $\Theta_r \Q[G]$. Hence Conjecture~\ref{Stark_over_Q} is equivalent to the statement that the isomorphism~(\ref{omega_log_isom}) maps $\E_r \Q \bigwedge^r U$ to $\Theta_r \Q[G]$. If the conjecture is true then $\varepsilon_r \in \E_r \Q \bigwedge^r U$. Suppose conversely that $\varepsilon_r \in \E_r \Q \bigwedge^r U$. Then $(\Omega \circ \log^{(r)})^{-1} (\Theta_r \Q[G]) = \Q[G] \varepsilon_r \subset \E_r \Q \bigwedge^r U$. In fact $\Q[G] \varepsilon_r = \E_r \Q \bigwedge^r U$ since these $\Q$-vector spaces have the same dimension $\# C_r$. Thus the conjecture is true.
\end{proof}

We consider this element $\varepsilon_r$ (conjecturally belonging to $\E_r \Q \bigwedge^r U$) as a generalization of the ``Stark unit''. Similarly to the latter, $\varepsilon_r$ is defined via units modulo torsion and it is unique as soon as some place of the bigger field $K$ over each place of the base field $k$ is fixed. 

There is a (certainly non-unique) way to write $c(\chi)$ as a single determinant. Let us show it for $\chi \ne {\bf 1}$. Consider $\E = \sum_{\alpha \in Gal(\Q(\chi)/\Q)} \E_{\chi^{\alpha}}$ and let $\{ i_1 < \dots < i_r \} = \{i : \chi|_{G_{v_i}} \equiv 1 \}$.
Since $\E \Q \bigwedge^r X$ is generated over $\Q[G]$ by the element $\E w_{i_1} \wedge \dots \wedge \E w_{i_r}$, each element $x \in \E \Q \bigwedge^r X$ can be represented (non-uniquely) in the form $x=x_1 \wedge \dots \wedge x_r$ with some $x_i \in \E \Q X$. Since $\Q U \cong \Q X$ as $G$-modules, we can represent
\[
\E \varepsilon_r = \epsilon_1 \wedge \dots \wedge \epsilon_r
\]
with some (non-unique) $\epsilon_i \in \E \Q U$. Then
\[
c(\chi) = \frac 1{\prod_{j=1}^r \#G_{v_{i_j}}} \det \left( \underset{g \in G}\sum \chi(g) \log|g \epsilon_k|_{w_{i_j}} \right)_{1\le k,j \le r}.
\]

\section{Example with the Hilbert class field of $\Q(\sqrt{229})$.}

The field $k=\Q(\sqrt{229})$ has class number $h=3$. To construct its Hilbert class field together with some units in it let us consider the polynomial
\begin{equation}\label{unit}
x^3-4 x+1.
\end{equation}
Let $K$ be its splitting field. It can have degree 3 or 6 over $\Q$. Since this equation has discriminant $d=229$, the splitting field should contain $k$. Hence $K$ is of degree 6 with $Gal(K /\Q) \cong S_3$, and $[K:k]=3$. We claim that $K$ is the Hilbert class field of $k$. Let us show $K/k$ is an unramified extension. Let $F=\Q(x^3-4 x+1)$. This qubic field is ramified at 229 only, so $K$ is also ramified at 229 only since it is the compositum of $k$ and $F$. Further, 
\[
x^3-4 x+1 = (x-29)^2 (x-171) \mod 229
\]
implies that the decomposition of $(229)$ into primes in $O_F$ has form $\frak p^2 \frak b$. Then $\frak p O_K = \frak P$ should be prime and $\frak b O_K = \frak B^2$ for a prime $\frak B$. This is because $K$ is a Galois field, so in the decomposition of $(229)=\frak P^2 \frak B^2$ all primes should have the same power. This power is the ramification index  $e_K(229)=2$. We have $e_K(229)=e_k(229)=2$, therefore $K/k$ is unramified.

Let $U=O_K^{\times}$ and $x \in U$ be the root of~(\ref{unit}). Notice that $x-2$ is also a unit then. Indeed, $(x-2)^3+6(x-2)^2+8(x-2)+1=x^3-4 x+1=0$. We could also take some unit $\epsilon \in O_k^{\times}$, and it is easy to check that $\epsilon$ together with the conjugates of  $x$ and $x-2$ under $G=Gal(K/k)$ generate the subgroup of finite index in $U$. Let $\E_1 = \frac 13 \sum_{g \in G} g$ and $\E_2 = 1-\E_1$. Then $\epsilon \in \E_1 U$ and $x, x-2 \in \E_2 U$. Thus the element 
\[
x \wedge (x-2)
\]
generates $\E_2 \Q \bigwedge^2 U$ over $\E_2 \Q[G]$.

We also need the Stickelberger element $\Theta_2$. We compute it using the Kronecker limit formula for real quadratic fields invented by Hecke. Expansion at $s=0$ of the partial zeta function $\zeta(g, s)$ (for $g \in G$) starts as
\[
\zeta(g, s) = - \frac{\log(\epsilon)}2 s + \rho(g) s^2 + o(s^2), \;\;\; s \To 0 
\]
where $\epsilon=15.066372975210...$ is the fundamental unit of $k$ and the coefficient $\rho(g)$ can be computed as follows. We take a quadratic form $Q$ of discriminant $229$ which represents the class of ideals mapped to $g \in G$ by the Artin map from the class field theory. Let $x > x'$ be its roots. Then
\[
\rho(g) = -\frac 14 \int_{-\log(\epsilon)}^{\log(\epsilon)} \log \left( y(v) |\eta(z(v))|^4 \right) dv + \mathrm{const}
\]
where $z(v)=\frac{x e^{\frac v 2}- i x' e^{-\frac v2}}{e^{\frac v2}-i e^{-\frac v2}}$, $y(v)$ is imaginary part of $z(v)$ and the constant is independent of $g$. This follows from the functional equation and formula for expansion near $s=1$ given with similar integrals in~\cite{Siegel} (page 90) or~\cite{Zagier} (page 162). The unity $1 \in G$ can be represented by the form $Q_0=x^2-15x-1$, and the other elements (each of which generates $G$) by $Q_1= 3 x^2-11 x-9$ and $Q_2=9 x^2-11 x+3$. The integral $\int_{-\log(\epsilon)}^{\log(\epsilon)} \log(y |\eta(z)|^4) dv$ equals $-16.230647798060277...$ for $Q_0$. Since $\zeta(g,s)=\zeta(g^{-1},s)$, the value of this integral for both $Q_1$ and $Q_2$ is the same number $-6.808829434682019...$ Therefore we find 
\[
\Theta_2 = \E_2 \sum_{g \in G} \rho(g) g = 1.570303060563043... \cdot (1-g/2-g^2/2)
\]
where $g \in G$ is any of the two generators in the last expression.

Let us enumerate the roots of~(\ref{unit}) as $x_0 = -2.1149075414...$, $x_1 = 0.2541016883...$, $x_2 = 1.8608058531...$ Let $w_{1,2}: K \To \R$ be the embeddings in which $x$ takes the value $x_0$ and $\sqrt{229}$ is positive and negative correspondingly. Let $g$ be a generator of $G$ such that $x$ take the value $x_1$ under $g w_1$. Since $Gal(K/\Q)$ acts as the permutation group $S_3$ on the set $(x_0, x_1, x_2)$ and $G$ is the subgroup of even permutations (cyclic shifts in this case), we see that $x$ will take the values $x_2$, $x_2$ and $x_1$ under $g^2 w_1$, $g w_2$ and $g^2 w_2$ correspondingly. Hence
\[
\log^{(2)} \left( x \wedge (x-2) \right) = \det \begin{pmatrix} \underset{i=0} \sum^2 \log|x_i|\;g^i & \underset{i=0} \sum^2 \log|x_i|\; g^{-i} \\ \underset{i=0} \sum^2 \log|x_i+2|\; g^i & \underset{i=0} \sum^2 \log|x_i+2|\; g^{-i} \\ \end{pmatrix} w_1 \wedge w_2
\]
\[
= 7.066363772533693... \cdot (g-g^2) \; w_1 \wedge w_2
\]
Let $\Omega = w_1^* \wedge w_2^*$ and
\[
\Omega \left(\log^{(2)} \left( x \wedge (x-2) \right)\right) = 7.066363772533693... \cdot (g-g^2). 
\]
Thus $\Theta_2 = - 1.570303060563043... \cdot \frac{(g - g^2)^2}2 =  \frac 19 (g^2-g) \; \Omega \left(\log^{(2)}\left( x \wedge (x-2) \right)\right)$ and
\[
\varepsilon_2 = \frac 19 (g^2-g) \; x \wedge (x-2) = \frac 19 \left(\frac{x^{g^2}}{x^g}\right) \bigwedge (x-2). 
\]

\end{document}